# DEHN SURGERY ON ARBORESCENT KNOTS AND LINKS – A SURVEY

Ying-Qing Wu

ABSTRACT. This article is solicited by C. Adams for a special issue of *Chaos, Solitons and Fractals* devoted to knot theory and its applications. We present some recent results about Dehn surgeries on arborescent knots and links.

In this survey we will present some recent results about Dehn surgeries on arborescent knots and links. Arborescent links are also known as algebraic links [Co, BoS]. The set of arborescent knots and links is a large class, including all 2-bridge links and Montesinos links. They have been studied by many people, see for example [Ga2, BoS, Mo, Oe, HT, HO].

We will give some definitions below. One is referred to [He] and [Ja] for more detailed background material for 3-manifold topology, to [Co, BoS, Ga2, Wu3] for arborescent tangles and links, to [Th1] for hyperbolic manifolds, and to [GO] for essential laminations and branched surfaces.

**0.1. Surfaces and 3-manifolds.** All surfaces and 3-manifolds are assumed orientable and compact, and surfaces in 3-manifolds are assumed properly embedded. Recalled that a surface $F$ in a 3-manifold $M$ is *compressible* if there is a loop $C$ on $F$ which does not bound a disk in $F$, but bounds one in $M$; otherwise $F$ is *incompressible*. A sphere $S$ in $M$ is a reducing sphere if it does not bound a 3-ball in $M$, in which case $M$ is said to be *reducible*. A 3-manifold is a *Haken* manifold if it is irreducible and contains an incompressible surface. $M$ is *hyperbolic* if it admits a complete hyperbolic metric. $M$ is *Seifert fibered* if it is a union of disjoint circles. Seifert fibered spaces are well understood, hyperbolic manifolds have also been deeply studied, although many aspects of which remain mysterious. Thurston's Geometrization Conjecture asserts that any closed irreducible manifold can be decomposed along incompressible tori into hyperbolic and Seifert fibered pieces, which has been proved by Thurston for Haken manifolds [Th2].

**0.2. Laminations.** Although Haken manifolds have many desirable properties, there do exist many non-Haken irreducible manifolds, which contain no incompressible surfaces. In [GO] Gabai and Oertel defined essential laminations. An

1991 *Mathematics Subject Classification.* 57M25, 57M50, 57N10.
Research at MSRI supported in part by NSF grant #DMS 9022140.

Typeset by $\mathcal{A}\mathcal{M}\mathcal{S}$-TEX

essential lamination in $M$ consists of noncompact incompressible surfaces satisfying certain technical conditions. A manifold containing an essential lamination is said to be *laminar*. Haken manifolds are also laminar. Gabai and Oertel showed that a laminar manifold has some nice properties of Haken manifolds, for example, it is irreducible, and its universal covering is $\mathbb{R}^3$. Laminations also play some other important rules. For example, the result of Brittenham [Br] shows that if the complementary regions of an essential lamination is not a product, then the manifold is not a small Seifert fibered space. This is a very important geometric criteria, and it played a crucial rule in the proof of the classification theorem of Dehn surgeries on 2-bridge knots, see Theorem 1.3 below.

**0.3. Dehn surgery.** Given a knot $K$ in the 3-sphere $S^3$, denote by $E(K)$ the manifold obtained by removing its regular neighborhood $N(K)$ from $S^3$, called the exterior of $K$. Its boundary is a torus $T$. Let $\gamma$ be an essential loop on $T$. Then the Dehn surgery manifold $K(\gamma)$ is obtained by gluing a solid torus $J$ back to $E(K)$ so that the meridian of $J$ is glued to the curve $\gamma$. Surgery with slope the meridian of $K$ will get $S^3$ back, and is called the trivial surgery. A surgery on a link $L$ is defined similarly. For our purpose, we define a nontrivial surgery on a link to be one which is nontrivial on all components of $L$. *All surgeries below are assumed nontrivial.* If $(m,l)$ is the standard meridian-longitude pair of $K$, then each slope $\gamma$ can be written as $pm+ql$ with $p,q$ coprime, which corresponds to a rational number $p/q$, or $\infty$ if $\gamma$ is the meridional slope (i.e. the trivial slope), hence nontrivial slopes for knots in $S^3$ are also denoted by rational numbers, see Rolfsen's book [Rf] for more details. It is known that all closed orientable 3-manifolds can be obtained by surgery on links in $S^3$.

**0.4. The main problems.** A surgery is a called *Haken, hyperbolic*, or *laminar surgery* if the resulting manifold is Haken, hyperbolic, or laminar, respectively. It is an *exceptional surgery* if the resulting manifold is reducible, toroidal, or is a small Seifert fibered space. A non exceptional surgery is also called a *hyperbolike* surgery [Gor], as the Geometrization Conjecture asserts that such a manifold is a hyperbolic manifold.

**Problem 0.1.** *Given a knot $K$ in $S^3$, how many surgeries on $K$ are non-Haken, non-laminar, or non-hyperbolic? How many are exceptional?*

There have been many results on these problems for surgery on knots. See [Gor], [Lu] and [Ga3] for excellent surveys and frontier problems. There is also a closely related problem asking how many Dehn fillings on a 3-manifold are nonsimple, see [Wu4] for some recent results and current status in this aspect. Recently Boyer and Zhang [BZ1,BZ2] have applied representation variety theory to make some important new progress. See Boyer's survey article in this volume for more details. For most arborescent knots, we have much better understanding, due to the special structures of these knots. This is also the case for links. There are very few results about Dehn surgery on general links concerning the above problems. However, we can show that surgeries on a large class of arborescent links always produce Haken manifolds.

**0.5. Arborescent tangles.** A *tangle* $(B, T)$, or simply $T$, is a pair with $B$ a 3-ball, and $T$ a properly embedded 1-manifold. We will only consider 2-string tangles, i.e. $T$ is a union of two arcs and possibly some circles. A *rational tangle* $T(r)$ corresponding to a rational number $r$ is obtained by first drawing two strings of slope $r$ on the boundary of a "pillowcase" $B$, then push the interior of the string into the interior of $B$. A *Montesinos tangle* $T(r_1, \ldots, r_n)$ is obtained by putting rational tangles $T(r_1), \ldots, T(r_n)$ together in a horizontal way. For example, the tangle in Figure 1(a) is the Montesinos tangle $T(1/2, 1/3, -1/3)$. We always assume that $r_i$ are non integers.

Figure 1

Given two tangles $(B_1, T_1)$ and $(B_2, T_2)$, we can make a new tangle $(B, T)$ by choosing a disk $D_i$ on $\partial B_i$ which intersects two ends of $T_i$, and gluing $D_1$ to $D_2$. Thus a Montesinos tangle is simply a sum of several rational tangles. In general, the sum of two tangles depends on the choice of the gluing disks and the gluing map, so the sum of two Montesinos tangles may no longer be a Montesinos tangle. An *arborescent tangle* is one that can be obtained by summing several Montesinos tangles together in an arbitrary order. For example, the tangle in Figure 1(b) is an arborescent tangle, which is the sum of two Montesinos tangles $M(1/2, 1/3)$ and $M(1/3, -1/3)$.

**0.6. Arborescent links.** Given a tangle $(B, T)$, we can add two arcs on the boundary of $B$, to form a knot or link in $S^3$. If $T$ is an arborescent tangle, then the corresponding link $L$ is called an *arborescent link*. The minimal number of rational tangles used to build $L$ is called the *length* of $L$. Another way to get an arborescent link is to take two arborescent tangles and glue the entire boundaries of the 3-balls to each other. One can see that these two definitions are equivalent. If $T$ is a rational tangle $T(r)$, we add two vertical lines on the boundary of the 3-ball (considered as a pillowcase), denote the resulting knot or link by $K_r$, and call it a *2-bridge knot or link*. If $T(r_1, \ldots, r_n)$ is a Montesinos tangle, we add two horizontal lines, and call the resulting link a *Montesinos link,* denoted by $M(r_1, \ldots, r_n)$. When $n = 1$, we have the relation $M(-1/r) = K_r$. It is known that $M(r_1, r_2)$ is always equivalent to some $K_r$, so we do not need to discuss Montesinos links of length 1 or 2.

For our purpose, we divide all arborescent knots into three types. Type I knots are 2-bridge knots or Montesinos knots of length 3. $K$ is a type II knot if it is the union of two arborescent tangles, each of which is a Montesinos tangle $M(1/2, r_i)$. All other arborescent knots are of type III. Similarly for links.

1. SURGERY ON 2-BRIDGE KNOTS AND LINKS

2-bridge knots are considered the simplest knots, and have been deeply studied. Incompressible surfaces in 2-bridge knot exteriors are classified by Hatcher and Thurston [HT], Property P is proved for 2-bridge knots by Takahashi [Ta]. The links $K_r$ with $r = 1/n$ are torus knots or links, whose complements are Seifert fiber spaces. All the other 2-bridge knots and links are hyperbolic, meaning that its complement admits a complete hyperbolic structure.

The work of Hatcher and Thurston [HT] shows that there is no closed incompressible surfaces in the exterior $E(K)$ of a 2-bridge knot $K$, except those parallel to the boundary of $E(K)$. By Hatcher's theorem [Ha], this implies that all but finitely many surgery manifolds are non-Haken. Moreover, all Haken surgery slopes can be calculated for any given 2-bridge knot. This settles the Haken surgery problem for 2-bridge knots. For laminar surgeries we have the following result of Delman. Recall that an essential lamination in the exterior of a knot $K$ is said to be *persistent* if it remains essential after all surgeries.

**Theorem 1.1.** [De1] *There exists a persistent lamination in the exterior of any nontorus 2-bridge knot.*

Thus surgeries on hyperbolic 2-bridge knots are always laminar.

To determine exceptional surgeries on 2-bridge knots, we need the following notations. A rational number $p/q$ associated to a 2-bridge knot $K = K_{p/q}$ has a partial fraction decomposition $p/q = [b_1, \ldots, b_n] = 1/(b_1 - 1/(b_2 - \ldots - 1/b_n)\ldots)$. We can always assume that $|b_i| > 1$. $K$ is a torus knot if and only if $n = 1$. The knot can be drawn as in Figure 2.

Figure 2

Surgeries on torus knots are well understood. They are all Seifert fibered except one, which is reducible. The simplest hyperbolic 2-bridge knot is the Figure 8 knot $K_{[2,-2]}$. Surgeries on this knot has been completely clarified by Thurston.

**Theorem 1.2.** [Th1] *Let $K$ be the Figure 8 knot. Then $K(\gamma)$ is hyperbolic for all but nine $\gamma$: $K(\gamma)$ is toroidal for $\gamma = 0, 4, -4$, and is Seifert fibered for $\gamma = -1, -2, -3, 1, 2, 3$.*

The following theorem of Brittenham and Wu classifies exceptional Dehn surgeries on all the other 2-bridge knots. It determines all toroidal or small Seifert fibered surgeries. There is no reducible surgery [HT, De1].

**Theorem 1.3.** [BW] *Let $K$ be a 2-bridge knot.*
*(1) If $K = K_r$ admits an exceptional surgery, then $r = [b_1, b_2]$ for some $b_1, b_2$.*
*(2) If $K = K_{[b_1,b_2]}$ with $|b_1|, |b_2| > 2$, then $K(\gamma)$ is exceptional for exactly one $\gamma$, which yields a toroidal manifold. When both $b_1$ and $b_2$ are even, $\gamma = 0$. If $b_1$ is odd and $b_2$ is even, $\gamma = 2b_2$.*
*(3) If $K = K_{[2n,\pm 2]}$ and $|n| > 1$, $K(\gamma)$ is exceptional for exactly five $\gamma$: $K(\gamma)$ is toroidal for $\gamma = 0, \mp 4$, and is small Seifert fibered for $\gamma = \mp 1, \mp 2, \mp 3$.*

A 2-bridge knot $K_r$ with $r = [2n, \pm 2]$ is called a *twist knot*. The theorem says that $K = K_{p/q}$ admits some exceptional surgeries if and only if it can be drawn as in Figure 2 with only two boxes. Such a knot admits either one or two toroidal surgeries, depending on whether it is a twist knot, and it admits some Seifert fibered surgeries if and only if it is a twist knot, which admits three such surgeries. The Figure 8 knot and the trefoil knot are also twist knots, but they have $|n| = 1$, and have been excluded from the theorem.

The main part of the proof of Theorem 1.3 is to show that certain surgered manifolds are not small Seifert fiber spaces. It uses the techniques of Delman [De1, De2] and Roberts [Rb] to construct essential laminations, which have the property that its complementary regions are not product. One can then apply Brittenham's criteria [Br] to show that the manifolds are not small Seifert fibered spaces.

**Problem 1.4.** *Show that the non exceptional surgeries for 2-bridge knots are really hyperbolic.*

Of course this would follow from the Geometrization Conjecture, or from the Orbifold Conjecture (see Kirby's problem set [Ki].) Using the techniques of Thurston [Th1], and Neumann and Reid [NR], it is possible to solve the problem for most twist knots.

The above result also gives some evidence to the following conjecture.

**Conjecture 1.5.** *If $K(\gamma_1)$ and $K(\gamma_2)$ are two integral nonhyperbolic surgeries, then $K(\gamma)$ is nonhyperbolic for all $\gamma$ between $\gamma_1$ and $\gamma_2$.*

We now consider 2-bridge links $L = K_{p/q}$. It is easy to see that $L$ has two components if and only if $q$ is even. There is no closed essential surfaces in the complement of such a link, so most surgeries are non-Haken [Fl]. Recall that we only consider nontrivial surgeries $L(\gamma_1, \ldots, \gamma_n)$, which means that all $\gamma_i$ are nontrivial slopes. The following theorem is a generalization of Delman's result (Theorem 1.2).

**Theorem 1.6.** [Wu2] *A 2-bridge link $L = K_{p/q}$ admits a non-laminar surgery if and only if $p/q = [m, n]$ for some odd integers $m$ and $n$.*

The proof shows that if $p/q \neq [m,n]$ then there is a persistent essential lamination in the link exterior. It is also known that a surgery on a hyperbolic 2-bridge link is laminar if both surgery slopes are non integral.

**Problem 1.7.** *Determine all non laminar surgeries and all exceptional surgeries on 2-bridge links of two components.*

There are some partial results. For example, the proofs of [BW] and [Wu4, Theorem 2.5] can be modified to show that if the 2-bridge link is "complicated" enough then it admits no exceptional surgeries. Also, if $K_{p/q}(n_1/m_1, n_2/m_2)$ is exceptional, then either some $m_i = 1$, or $p/q = [r, s]$ and some $m_i = 2$. Given such a link, one can use Jeff Weeks' computer program *Snappea* [We] and Casson's program *geo* [Ca] to get some sense about which surgeries on it are exceptional.

## 2. Surgery on Montesinos knots of length 3

Let $K = M(r_1, r_2, r_3)$ be a Montesinos knot of length 3. There is no closed incompressible surface in the exterior of such a knot, hence most surgeries are non-Haken. Delman has constructed persistent essential laminations in the complements of most such knots.

**Theorem 2.1.** [De2] *The exterior of $K$ contains a persistent essential lamination unless $K = M(x, 1/p, 1/q)$ or its mirror image, where $x \in \{-1/2n, -1 \pm 1/2n, -2 + 1/2n\}$, and $p, q, n$ are positive integers.*

The $(-2, 3, 7)$ pretzel knot is a well known example. In the above notation it is the knot $M(-1/2, 1/3, 1/7)$. It admits 6 exceptional surgeries, see [BH] for more details. It is not known if it has any other exceptional surgeries, although *Snappea* suggests that those are all. Some other knots in Theorem 2.1 also admit exceptional surgeries, but it is possible that for most of the knots listed in the theorem, all surgeries are still laminar and non exceptional. It is a very interesting problem to construct essential laminations for Dehn surgeries on these remaining knots.

Using the results of [Br] and the techniques of [Wu4] it can be shown that if $\gamma$ is a nonintegral slope and if $K$ is a knot admitting a persistent lamination given by Theorem 2.1, then $K(\gamma)$ is atoroidal, and is not a small Seifert fibered space; hence all exceptional surgeries on $K$ are integral surgeries. Also, if $K$ is complicated enough so that the essential lamination constructed by Delman has 4 cusps around the knot, then the above is true for integral slopes as well, and $K$ admits no exceptional surgeries.

Thus, unlike the other types of arborescent knots, our knowledge about surgery on Montesinos knots of length 3 are less satisfactory. The laminar surgery problem is open for those listed in Theorem 2.1, and the exceptional surgery problem is open for even more knots. Surgery on length 3 Montesinos links has never been touched. There are still a lot of things to be worked out about these knots and links.

## 3. Surgery on type II and type III knots

Recall that a link is of type II if it is the union of two arborescent tangles of the form $T(1/2, r_i)$, and all the other arborescent links are either a Montesinos links of length at most 3 or type III links. Surgery on type II and type III knots are much better understood, due to the fact that there are some closed incompressible surfaces in the knot exterior. Earlier, Oertel [Oe] showed that surgeries on Montesinos knots of length $\geq 4$ are Haken. The following theorem gives a satisfactory solution to our problems for all surgeries on type III knots and all nonintegral surgeries on type II knots.

**Theorem 3.1.** [Wu1] *If $K$ is a type III arborescent knot, or if $K$ is of type II and $\gamma$ is a non-integral slope, then $K(\gamma)$ is a hyperbolic Haken manifold.*

In particular, all surgeries on a Montesinos knot $K = M(r_1, \ldots, r_n)$ with $n \geq 4$ are Haken and hyperbolic because they are of type III. (A type II Montesinos link has at least 2 components.)

The proof that $K(\gamma)$ is hyperbolic uses Thurston's geometrization theorem for Haken manifolds. It was shown that the surgered manifolds are Haken, atoroidal, and are not small Seifert fibered manifolds. Theorem 3.1 is not true for integral surgeries on type II knots. Each type II knot exterior contains infinitely many isotopy classes of connected, closed, incompressible surfaces, but none of them can survive under any integral surgery.

**Theorem 3.2.** [Wu1] *If $K$ is a type II arborescent knot, then all closed incompressible surfaces in $E(K)$ are compressible in $K(\gamma)$ for any integral slope $\gamma$.*

**Problem 3.3.** *Are there any (integral) exceptional surgery on type II knots?*

Haken manifolds are laminar because an incompressible surface can be considered as an essential lamination, so those surgeries in Theorem 3.1 are automatically laminar. It turns out that integral surgeries on type II knots are also laminar, which gives a complete solution to the laminar surgery problem for these knots.

**Theorem 3.4.** [Wu1] *If $K$ is an arborescent knot of type II or III, then $K(\gamma)$ is laminar for all slopes $\gamma$.*

**Corollary 3.5.** [Wu1] *The Property P conjecture and the Cabling Conjecture are true for all arborescent knots.*

These conjectures are still among the most interesting problems in Dehn surgery theory. The *property P conjecture* asserts that surgery on a nontrivial knot will not produce a homotopy sphere. Gordon and Luecke's Knot Complement Theorem [GL] implies that this would follow from the Poincare conjecture. The conjecture has also been proved for satellite knots [Ga1], and symmetric knots [CGLS]. Recently Delman and Roberts [DR] proved it for alternating knots. The *Cabling Conjecture* asserts that all surgeries on hyperbolic knots are irreducible. It has been proved for satellite knots [Sch], alternating knots [MT], strongly invertible knots [Eu], and those with bridge number at most 4 [GL2].

The above theorems are proved using tangles. They are consequences of the following general result: *If $K$ is the union of two nontrivial atoroidal tangles, then $K(\gamma)$ are laminar for all $\gamma$. Moreover, if at least one of the tangle space has incompressible boundary, then $K(\gamma)$ are also hyperbolic and Haken.* In particular, this applies to all knots which have nontrivial arborescent parts [BoS].

Some results of [Wu3] are used in the proof. A tangle is said to be *hyperbolic* if its exterior admits a hyperbolic structure with totally geodesic boundary. This is equivalent to saying that there is no essential surface of nonnegative Euler number in the tangle space $E(T)$, so $T$ is also called a *simple tangle*. All hyperbolic arborescent tangles with no loops have been classified [Wu3]. The following is for the easy case when $T$ is a Montesinos tangle.

**Theorem 3.6.** [Wu3] *Let $M(r_1, \ldots, r_n)$ be a Montesinos tangle containing no loops. Suppose $r_i = p_i/q_i$ with $|p_i| \leq q_i/2$.*

*(a) If $n = 2$, $M(r_1, r_2)$ is nonhyperbolic if and only if either both $p_i = \pm 1$, or one of the $r_i = 1/2m$ for some integer $m$.*

*(b) If $n \geq 3$, $M(r_1, \ldots, r_n)$ is nonhyperbolic if and only if either $r_1$ or $r_n$ is of the form $1/2m$ with $|m| \geq 2$.*

Note that we can always arrange by flyping moves so that $|p_i| \leq q_i/2$. It is very easy to use this theorem to determine if a Montesinos tangle is hyperbolic, which can also be done for any arborescent tangle using the classification theorem of [Wu3]. However, if $T$ is not an arborescent tangle, then it may be very hard to decide if it is hyperbolic. It is not too hard to show that the two tangles in Figure 3(a,b) are nonhyperbolic, but it seems hard to work it out for those in Figure 3(c,d).

**Problem 3.7.** *(1) Generalize Theorem 3.6 and the classification theorem of [Wu3] to arborescent tangles containing loops.*

*(2) Find a practical algorithm to decide if a given tangle is hyperbolic.*

Figure 3

## 4. Surgery on arborescent links of multiple components

Most research on surgery has been done only for knots. Usually the results are not readily generalized to surgery on links of multiple components, a major difficulty being that surgery on one component of the link may change the property of the others. There are some exceptions: One is Thurston's hyperbolic surgery theorem [Th1], which says that if $L$ is a hyperbolic link, then when excluding finitely many slopes on each component of $L$, all the remaining surgeries are hyperbolic. For exceptional surgeries, the $2\pi$-theorem of Gromov-Thurston says that it suffices to exclude 48 slopes from each component, see [BH]. Another interesting result is Scharlemann's simultaneous crossing change theorem, see [Sch2].

We study the problem of which arborescent link $L$ has the property that all nontrivial surgeries are Haken. We must assume that $L$ is not a Montesinos link of length at most 3, because otherwise there is no closed incompressible surface in the link exterior, so most surgeries are non-Haken. Since there is no limit on the number of components in $L$, the situation is more complicated. Examples show that "bad" surgery may occur if some of the rational tangles are $1/2$ tangles, even if the length of $L$ is arbitrarily large. Also, the Montesinos link $M(1/2k_1, \ldots, 1/2k_n)$ admits a reducible surgery. These examples show that the conditions in the following theorem cannot be removed.

**Theorem 4.1.** [Wu2] *Let $L = k_1 \cup \ldots \cup k_n$ be an arborescent link of type II or III, so that (1) $L$ is composed of rational tangles $T(p_i/q_i)$ with $q_i > 2$, and (2) no length 2 tangle of $L$ is a Montesinos tangle of the form $T(1/2q_1, 1/2q_2)$. Then all nontrivial surgeries on $L$ produce Haken manifolds.*

The proof of this theorem involves essential branched surfaces. The idea is to cut $S^3$ into pieces, each containing at most one component of the link. Thus the branched surface plays the rule of an insulator separating the components of $L$, which reduces the problem to surgery on knots in 3-manifolds.

When $L$ has only two components, the second condition in Theorem 4.1 can be removed.

**Theorem 4.2.** [Wu2] *Suppose $L = l_1 \cup l_2$ is a two component arborescent link of type II or III. If none of the rational tangles used to build $L$ is a $(1/2)$-tangle, then all nontrivial surgeries on $L$ produce Haken manifolds.*

It is likely that all the surgeries in the above theorems are hyperbolic. Certainly they are not small Seifert fibered, but we don't know if any of them are toroidal.

**Problem 4.3.** *Determine if all the surgeries in the above theorems are atoroidal.*

The *Strong Property P Conjecture* of Gabai [Ga3] asserts that all knots in $S^3$ have Strong Property P, i.e, $K(\gamma)$ does not contain fake 3-ball for all $K$ and $\gamma$. By Gordon and Luecke's Knot Complement Theorem [GL], this implies the Property P Conjecture. Since all closed 3-manifolds are obtained by surgery on links, the Poincare Conjecture is equivalent to the following link version of Strong Property P Conjecture.

**Conjecture 4.4.** *All links $L$ in $S^3$ have strong property P, that is, surgeries on $L$ produce manifolds containing no fake 3-balls.*

To prove the Poincare conjecture, it suffices to prove Strong Property P for a class of links which is large enough so that all closed 3-manifolds can be obtained by surgery on a link in the class. The following theorem shows that the conjecture is true for Montesinos links.

**Theorem 4.5.** *All Montesinos links have Strong Property P.*

*Proof.* A link $L$ is said to be *strongly invertible* if there is an involution $\varphi$ of $S^3$ with axis $C = \text{Fix}(\varphi)$ intersecting each component of $L$ twice, and $\varphi(L) = L$. It can be shown that if $L$ is strongly invertible, then $L(\gamma)$ is a double branched cover of $S^3$. In this case one can use the Smith Conjecture [MB] to show that $L(\gamma)$ contains no fake 3-ball for any $\gamma$. Thus all strongly invertible links have Strong Property P. The result now follows from Corollary 3.5 and Theorem 4.2, with some simple argument showing that a Montesinos link is either strongly invertible or is a link of at most two components satisfying the conditions of Theorem 4.2. □

**Problem 4.6.** *Prove Strong Property P for all arborescent links, or some other classes of links.*

Department of Mathematics, University of Iowa, Iowa City, IA 52242; and MSRI, 1000 Centennial Drive, Berkeley, CA 94720-5070
*E-mail address*: `wu@math.uiowa.edu, wu@msri.org`